\documentclass[12pt]{article}
\usepackage{amssymb}
\usepackage{amsthm}
\usepackage{amsmath}
\usepackage{graphicx}
 \newtheorem{thm}{Theorem}[section]
 
 \newtheorem{lem}[thm]{Lemma}
 
 \theoremstyle{definition}
 
 \newtheorem{rem}[thm]{Remark}
 \numberwithin{equation}{section}

\begin{document}
\title{On 2D Viscoelasticity with Small Strain}
\author{Zhen Lei
\thanks{School of Mathematical Sciences, Fudan University,
Shanghai 200433, P. R. China. {\it email: leizhn@yahoo.com,
zlei@fudan.edu.cn}}}
\date{}
\maketitle

\begin{abstract}
An exact two-dimensional rotation-strain model describing the
motion of Hookean incompressible viscoelastic materials is
constructed by the polar decomposition of the deformation tensor.
The global existence of classical solutions is proved under the
smallness assumptions only on the size of initial strain tensor.
The proof of global existence utilizes the weak dissipative
mechanism of motion, which is revealed by passing the partial
dissipation to the whole system.
\end{abstract}

\section{Introduction}

Viscoelastic materials include a wide range of fluids with elastic
properties, as well as solids with fluid properties. Let us start
with the standard description of general mechanical evolutions to
introduce some notations and definitions. Any deformation can be
represented by a flow map (particle trajectory) $x(t,X)$, $0 \leq
t < T$, a time-dependent family of orientation-preserving
diffeomorphisms. $X$ is the original labeling (Lagrangian
coordinate) of the particles. $X$ is also referred to as the
material/reference coordinate. $x$ is then the observer's
(Eulerian) coordinate. In general, the velocity field $u(t,x)$ is
defined as the time derivative of the flow map, which is a vector
field defined on the Eulerian coordinate:
\begin{equation}\nonumber
u(t,x) = x_t\big(t, X(t, x)\big).
\end{equation}
The deformation tensor $\widetilde{F}$ in Lagrangian coordinate is
defined by
\begin{equation}\nonumber
\widetilde{F}_{ij}(t, X) = \frac{\partial x_i}{\partial X_j}(t,
X).
\end{equation}
When one uses the Eulerian description, the deformation gradient
tensor $F$ in spatial coordinate is defined by
\begin{equation}\nonumber
F(t, x) = \widetilde{F}\big(t, X(t, x)\big).
\end{equation}

Recently, an extensively studied system modelling the motion of
isotropic incompressible viscoelastic materials of Oldroyd-B type
takes:
\begin{equation}\label{a3}
\begin{cases}
\nabla \cdot u = 0,\\
  u_t + u \cdot \nabla u + \nabla p
      = \mu\Delta u + \nabla\cdot(\frac{\partial W(F)}{\partial
      F}F^T),\\[-4mm]\\
  F_t + u \cdot \nabla F = \nabla uF,
\end{cases}
\end{equation}
where $p(t, x)$ is the pressure, $\mu$ is the viscosity and $W(F)$
is the elastic energy function. The global existence of classical
solutions near equilibrium to \eqref{a3} in the 2D Hookean case
was established by Lin, Liu and Zhang \cite{Lin2} (see also Lei
and Zhou \cite{Lei5} via the incompressible limit method). Later
on, Lei, Liu and Zhou \cite{Lei4} proved the global existence of
classical small solutions to system \eqref{a3} in both 2D and 3D
cases (see also Chen and Zhang \cite{ChenZhang2006}). Very
recently, an improved result was obtained by Kessenich \cite{Kess}
by removing the dependence of the smallness of the initial data on
the viscosity via the hyperbolic energy method.

To study the different contributions of strain and rotation parts,
Friedrich \cite{Friedrichs} observed that the smallness of the
strain in nonlinear elasticity can be realized through the polar
decomposition of the deformation tensor. John \cite{John1, John2}
showed that no pointwise estimate for rotations in terms of
strains can exist even in the case of small strain. In the work of
Friedrich and John, no PDE is involved. From the PDE point of
view, Liu and Walkington \cite{Liu1} considered approximating
systems resulting from the special linearization of the original
system with respect to the strain. Lei, Liu and Zhou \cite{Lei3}
constructed a 2D rotation-strain viscoelastic model and proved the
global existence their classical solutions with small strain. One
shortcoming of the model in \cite{Lei3} is that the dynamics of
strain and rotation are not equivalent to that of the deformation
tensor.

One main concern of this paper is to construct an exact
rotation-strain model for motions of Hookean viscoelastic
materials of Oldroyd-B type, which takes
\begin{equation}\label{b}
\begin{cases}
\nabla\cdot u = 0,\\
u_t + u\cdot\nabla u + \nabla p = \mu\Delta u + \nabla\cdot (VV^T)
      + 2\nabla\cdot V,\\[-4mm]\\
V_t + u\cdot\nabla V = D(u) + \frac{1}{2}(\nabla uV +
      V\nabla u^T)\\
\quad\quad+\ \frac{1}{2}\omega_{12}(u)(VA - AV)
      - \frac{1}{2}\gamma(VA - AV),\\[-4mm]\\
\nabla^\perp\theta_t + u\cdot\nabla\nabla^\perp\theta =
      \frac{1}{2}\Delta u - \begin{pmatrix}- \nabla_2u\cdot\nabla\theta\\
      \nabla_1u\cdot\nabla\theta\end{pmatrix} + \nabla^\perp\gamma.
\end{cases}
\end{equation}
Here $V$ is the strain part of the deformation tensor, $\theta$ is
the rotation angle and $\gamma$ is a nonlinear term. The
derivation of this model will be presented in section 2. It seems
that  the velocity field $u$ is the only dissipative unknown. We
will prove the global existence of classical solutions to
\eqref{b} for small initial strain $V$ and initial gradient of
$\theta$ ($\theta$ may be arbitrarily large). The proof utilizes
the weak dissipative mechanism of system \eqref{b} by passing the
dissipation of the velocity field $u$ to the strain tensor $V$ and
the gradient of the angle variable $\nabla\theta$. Note that the
equations for $u$ and $V$ have already formed an closed dynamical
system. However, our proof needs the use of the equation for
$\theta$ to carry out the weak dissipative mechanism of the whole
system. Without the help of the equation for $\theta$ and the
underlying physical relationship between $V$ and $\theta$ (see
Lemma \ref{lem2}), the corresponding analysis issue is still open.

At last, let us also cite some related results on elastodynamics
when the viscosity $\mu$ is zero. The global existence of 3D
compressible elastodynamics was established by Agemi \cite{Agemi}
and Sideris \cite{Sideris1} independently under the null condition
and small initial data. An earlier almost global existence result
can be founded in John and Klainerman \cite{John3} (see also
Klainerman and Sideris \cite{Klainerman3}). Hyperbolic energy
method involving Klainerman's vector field developed by Klainerman
\cite{Klainerman2} (see also Christodoulou \cite{Christodoulou})
plays an essential role in \cite{Klainerman3, Sideris1}. The
incompressible case was later on studied by Sideris and Thomases
in \cite{SiderisThomases2005, SiderisThomases2006,
SiderisThomases2007}. A natural question to ask is that what
happens for the solutions of \eqref{a3} when the viscosity
vanishes, which is answered very recently by Kessenich \cite{Kess}
in the case of 3D. Other results on Oldroyd-B model can be found
in \cite{Chemin, Lei1, Lin1, Lions2}.

The remaining part of this paper is organized as follows. In
section 2 we present the derivation of the rotation-strain model.
Then the main result of this paper is stated. Section 3 is devoted
to exploring some of the intrinsic properties of the viscoelastic
system. The proof of global existence theorem is completed in
section 4.

\section{Derivation of the Strain-Rotation Model}

Constitutive relations generally involve the deformation tensor
$F$. In particular, the polar decomposition of $F$ is important.
For any non-singular matrix $F$, it is easy to see that $FF^T$ is
positive definite symmetric. Thus, there exists a unique positive
definite symmetric matrix $I + V$ such that $(I + V)^2 = FF^T$,
where $I$ is the $2 \times 2$ identity matrix. Then the matrix $R
= (I + V)^{- 1}F$ is orthogonal. Consequently the non-singular
matrix $F$ can be decomposed uniquely in the form
\begin{equation}\label{b1}
F = (I + V)R,
\end{equation}
where $R$ is orthogonal
$$RR^T = I$$
and $I + V$ positive definite symmetric
\begin{equation}\label{b3}
V = V^T.
\end{equation}
Physically, this means that the deformation is decomposed into
stretching and rotation. Following a suggestion by K. O.
Friedrichs \cite{Friedrichs}, the tensor $I + V$ is called the
left stretch tensor, $V$ the strain matrix and $R$ the rotation
matrix.

Our first goal is to formulate the viscoelastic system in terms of
the velocity $u$ and the strain $V$. We first look at the momentum
equation, the second equation in \eqref{a3}. We restrict our
discussions below to a special elastic energy functional of the
Hookean linear elasticity $W(F) = \frac{1}{2}|F|^2$, it does not
reduce the essential difficulties for mathematical analysis.
Indeed, all the results we describe here can easily be generalized
to a more general class of elastic energy functionals of the
deformation tensor $F$ (see \cite{Lei4}). In terms of the strain
matrix $V$, the momentum equation takes the form of
\begin{equation}\label{b2}
u_t + u\cdot\nabla u + \nabla p = \mu\Delta u + \nabla\cdot (VV^T)
+ 2\nabla\cdot V
\end{equation}
where we used the equation \eqref{b1} and \eqref{b3}. We point out
here that through this paper, $\nabla$ with no indices stands for
the derivative with respect to the spatial coordinate $x$. For a
matrix-valued function $A$, $(\nabla\cdot B)_i$ means
$\nabla_jB_{ij}$, where summation over repeated indices will
always be well understood.

In order to get a dynamical equation for $V$, we plug \eqref{b1}
into the third equation of system \eqref{a3}. Denote the material
derivatives of a quantity $Q$ by
\begin{equation}\nonumber
\dot{Q} = Q_t + u\cdot\nabla Q.
\end{equation}
We have
\begin{equation}\nonumber
\dot{V}R + (I + V)\dot{R} = \nabla u(I + V)R.
\end{equation}
Noting that $R$ is orthogonal, we deduce that
\begin{equation}\label{b6}
\dot{V} + (I + V)\dot{R}R^T = \nabla u(I + V).
\end{equation}
By transposing the above equality \eqref{b6}, it is rather easy to
see
\begin{equation}\label{b7}
\dot{V} + R\dot{R}^T(I + V) = (I + V)\nabla u^T.
\end{equation}

Now let us write the orthogonal matrix $R$ as
\begin{equation}\nonumber
R = \begin{pmatrix}\cos\theta & - \sin\theta\\
    \sin\theta & \cos\theta\end{pmatrix}.
\end{equation}
It is immediately verified that
\begin{eqnarray}\nonumber
\dot{R}R^T = \begin{pmatrix}- \sin\theta & - \cos\theta\\
    \cos\theta & - \sin\theta\end{pmatrix}
             \begin{pmatrix}\cos\theta & \sin\theta\\
    - \sin\theta & \cos\theta\end{pmatrix}\dot{\theta}
= \begin{pmatrix}0 & - 1\\
        1 & 0\end{pmatrix}\dot{\theta}.
\end{eqnarray}
Define
\begin{equation}\label{b10}
A = \begin{pmatrix}0 & - 1\\
    1 & 0\end{pmatrix},
\end{equation}
we have
\begin{equation}\label{b11}
\dot{R}R^T = A\dot{\theta},\quad VA + AV = A{\rm tr}V \quad A^T =
A^{- 1} = - A.
\end{equation}
Thus, it follows from \eqref{b6}, \eqref{b7} and \eqref{b11} that
\begin{eqnarray}\nonumber
  (I + V)A\dot{\theta} - A^T(I + V)\dot{\theta} = \nabla u(I + V) - (I + V)\nabla u^T
\end{eqnarray}
which  is equivalent to
\begin{eqnarray}\nonumber
  \dot{\theta}I &=& - \frac{1}{2 + trV}A[\nabla u(I + V)
      - (I + V)\nabla u^T]\\\nonumber
  &=& - \frac{1}{2 + trV}[2\omega_{12}(u) + \nabla_ku_1V_{k2}
  - \nabla_ku_2V_{k1}]I\\\nonumber
  &=& - \omega_{12}(u)I + \gamma I,
\end{eqnarray}
where
\begin{eqnarray}\label{b14}
\omega(u) = \frac{1}{2}(\nabla u - \nabla u^T)
\end{eqnarray}
is the virticity tensor and
\begin{equation}\label{b02}
  \gamma = \frac{1}{2 + trV}[trV\omega_{12}(u)
      - (\nabla_ku_1V_{k2} - \nabla_ku_2V_{k1})].
\end{equation}
Consequently, we obtain that
\begin{equation}\label{b13}
\dot{\theta} = - \omega_{12}(u) + \gamma.
\end{equation}

Next, add up \eqref{b6} and \eqref{b7} and then insert \eqref{b11}
and \eqref{b13} into the resulting equation, we obtian
\begin{eqnarray}\label{b15}
\dot{V} &=& \frac{1}{2}[\nabla u(I + V) + (I + V)\nabla
                       u^T]\\\nonumber
        &&\quad -\ \frac{1}{2}[(I + V)\dot{R}R^T + R\dot{R}^T(I +
                       V)]\\\nonumber
        &=& \frac{1}{2}(\nabla u + \nabla u^T) + \frac{1}{2}(\nabla uV + V\nabla
           u^T)\\\nonumber
        &&\quad +\ \frac{1}{2}[(I + V)A + A^T(I + V)]\omega_{12}(u)\\\nonumber
        &&\quad -\ \frac{1}{2}\gamma[(I + V)A + A^T(I + V)]\\\nonumber
        &=& \frac{1}{2}(\nabla u + \nabla u^T) + \frac{1}{2}(\nabla uV + V\nabla
           u^T)\\\nonumber &&\quad +\ \frac{1}{2}\omega_{12}(u)(VA - AV) - \frac{1}{2}\gamma(VA - AV),
\end{eqnarray}
where in the last equality we used \eqref{b11}. If we denote the
symmetric part of the gradient of velocity by
\begin{eqnarray}\label{b16}
D(u) = \frac{1}{2}(\nabla u + \nabla u^T),
\end{eqnarray}
we can rewrite the above equation \eqref{b15} as
\begin{eqnarray}\label{b17}
 V_t &+& u\cdot\nabla V = D(u) + \frac{1}{2}(\nabla uV +
V\nabla u^T)\\\nonumber &+& \frac{1}{2}\omega_{12}(u)(VA - AV) -
\frac{1}{2}\gamma(VA - AV).
\end{eqnarray}

In order to formulate an incompressible dynamical problem, the
equations of mass conservation $\nabla\cdot u = 0$ and balance of
momentum \eqref{b2} supplemented by the transport equation
\eqref{b17} have in fact made up of a closure system. However, to
prove the global existence, we must include an equation for
$\nabla\theta$ to carry out some intrinsic properties of the
rotation-strain viscoelastic system, which will be seen more clear
in section 4.

Let
\begin{equation}\label{b19}
\nabla^\perp = \begin{pmatrix}- \nabla_2\\ \nabla_1\end{pmatrix}
\end{equation}
and then apply $\nabla^\perp$ to \eqref{b13} to yield
\begin{equation}\label{b20}
\nabla^\perp\theta_t + u\cdot\nabla\nabla^\perp\theta =
\frac{1}{2}\Delta u - \begin{pmatrix}- \nabla_2u\cdot\nabla\theta\\
\nabla_1u\cdot\nabla\theta\end{pmatrix} + \nabla^\perp\gamma.
\end{equation}
This equation will be used to explore the weak dissipation of
$\nabla\theta$, which eventually shows that $V$ is also weakly
dissipative. The equations \eqref{b2}, \eqref{b17} and \eqref{b20}
form the rotation-strain model \eqref{b} for viscoelastic fluids.
For future use, we also need an equation of $\nabla\theta$. For
this purpose, we apply the gradient operator $\nabla$ to
\eqref{b13} to get
\begin{equation}\label{b21}
\nabla\theta_t + u\cdot\nabla\nabla\theta =
- \begin{pmatrix}\nabla_1u\cdot\nabla\theta\\
\nabla_2u\cdot\nabla\theta\end{pmatrix} - \nabla\omega_{12}(u) +
\nabla\gamma.
\end{equation}

We will consider the Cauchy problem or the periodic boundary-value
problem of system \eqref{b}. The initial data takes of the form
\begin{equation}\label{bb}
u(0, x) = u_{0}(x),\ \ \ V(0, x) = V_{0}(x),\ \ \ \nabla\theta(0,
x) = \nabla\theta_0(x),\ \ \ x \in \Omega.
\end{equation}
where $\Omega \subseteq R^2$ is two dimensional torus or the
entire space. The above initial data will be imposed on the
following constraints
\begin{equation}\label{bbb}
\begin{cases}
   \nabla\cdot u_0 = 0,\\
   \det(I + V_0) = 1,\\
   \nabla\cdot V_0 = A(I + V_0)\nabla\theta_0.
\end{cases}
\end{equation}
The first two are just the consequences of the incompressibility
and the last one is understood as the consistency condition for
changing of variables (see Lemma \ref{lem1}, Lemma \ref{lem2},
Remark \ref{rem2} and Remark \ref{rem3} in the next section).

\begin{rem}\label{rem1}
In our earlier work \cite{Lei3}, the rotation matrix $R$ satisfies
a spontaneous but specified transport equation:
$$R_t + u\cdot\nabla R = - AR\omega_{12}(u) = \omega(u)R.$$
Consequently, it follows that
$$\theta_t + u\cdot\nabla\theta = - \omega_{12}.$$
This combining the third equation of system \eqref{a3} gives a
non-symmetric transport equation for $V$. While in the model
\eqref{b}, noting \eqref{b11} and with the aid of \eqref{b13}, it
is rather easy to find that
\begin{equation}\label{b18}
R_t + u\cdot\nabla R = AR\dot{\theta} = - AR\omega_{12}(u) +
\gamma AR.
\end{equation}
The last term of the right side of the above equation or the one
$\gamma$ in \eqref{b13} represents corrections of our previous
model, which leads to that the strain matrix must be symmetric.
However, the corrections reflect the intrinsic properties of
motions of viscoelastic fluids and the underlying physical
origins.
\end{rem}

For the rotation-strain model \eqref{b}, we will prove the
following theorem:

\begin{thm}\label{thm1}
Consider the Cauchy problem or the periodic initial-boundary value
problem for the rotation-strain viscoelastic model \eqref{b} and
\eqref{bb} with the intrinsic physical constraints \eqref{bbb} on
the initial data. Then there exists a unique global classical
solution $(u, V, \nabla\theta)$ which satisfies
\begin{eqnarray}
\nonumber&&\|u\|_{H^2(\Omega)}^{2} + \|V\|_{H^2(\Omega)}^{2} +
\|\Delta u + \frac{1}{\mu}\nabla^\perp\theta\|^2 +
\mu\int_0^\infty\Big[\|\nabla u\|_{H^2}^2\\\nonumber &&+\
\|\nabla(\Delta u + \frac{1}{\mu}\nabla^\perp\theta)\|^2 +
\|\nabla(\Delta u + \frac{1}{\mu}\nabla\cdot V)\|^2\Big]dt \leq
\frac{\mu^4(1 + \mu^4)}{C^2(1 + \mu^{10})},
\end{eqnarray}
if the initial data $u_0,\ V_0 \in H^2(\Omega)$ and
$\nabla\theta_0 \in H^{1}(\Omega)$ and
$$\|u_{0}\|_{H^{2}(\Omega)}^2 + \|V_{0}\|_{H^{2}(\Omega)}^2
+ \|\nabla\theta_{0}\|_{H^{1}(\Omega)}^2 < \frac{\mu^8}{M(1 +
\mu^{10})},$$ where $\Omega \subseteq R^2$ is two dimensional
torus or the entire space, $C$ and $M$ $(M > 2C^3)$ are big enough
constants independent of $t$ and $\mu$.
\end{thm}

The proof of the above theorem relies on a local existence theorem
and a \emph{priori} energy estimates. The proof of local existence
is standard and can be similarly done as in \cite{Lei4, Lin1},
thus is omitted here. Below we will only present the \textit{a
priori} estimate in Theorem \ref{thm1}.

\begin{rem}
In the case of $\Omega = \mathrm{T}^2$ being a two-dimensional
torus, we in fact recovered the results in our earlier work
\cite{Lei4}. To see this, we set
$$\overline{\theta}_0 = \frac{1}{|\mathrm{T}^2|}\int_{\mathrm{T}^2}\theta_0 dx.$$
By Poincar$\acute{e}$'s inequality, we have
\begin{equation}\label{b30}
\|\theta_0 - \overline{\theta}_0\|_{L^\infty} \leq C\|\theta_0 -
\overline{\theta}_0\|_{H^2} \leq C\|\nabla\theta_0\|_{H^1}.
\end{equation}
Thus, if we define
$$R_0 = \begin{pmatrix}\cos\overline{\theta}_0 & - \sin\overline{\theta}_0\\
\sin\overline{\theta}_0 & \cos\overline{\theta}_0 \end{pmatrix},$$
then
$$F(0, x) = \big[(I + V)R\big]_{t = 0} = R_0 + E_0$$
where
$$E_0 = \big[R - R_0 + VR\big]\big|_{t = 0}$$
is a small disturbance near a constant equilibrium $R_0$ by
\eqref{b30}, since $\|\nabla\theta_0\|_{H^1}$ is assumed to be
small in Theorem \ref{thm1}. In \cite{Lei4}, we studied the
general viscoelastic model including an equation for the
deformation tensor $F$ and proved the global existence of
classical solutions near constant equilibrium in both two and
three-dimensional spaces. However, our results and proofs here are
still of importance even in the periodic case because they provide
us a better understanding of the physical background of the
system.
\end{rem}


\section{Special Structures of the System}

In this section, we will explore some of the intrinsic properties
of the viscoelastic system mentioned in section 2. These
properties reflect the underlying physical origins of the problem
and in the meantime, and are also essential to the proof of the
global existence result here.

The following lemma will be used to explore the weak dissipation
of $V$.
\begin{lem}\label{lem1}
Assume that the second equality of \eqref{bbb} is satisfied and
$(u, V)$ is the solution of system \eqref{b}. Then the following
is always true:
\begin{equation}\label{c1}
\det (I + V) = 1
\end{equation}
for all the latter time $t \geq 0$. In other words, we have
\begin{equation}\label{coreq}
trV = - \det V.
\end{equation}
\end{lem}

$Proof.$ Using the identity that for any non-singular matrix $F$,
$$\frac{\partial\det F}{\partial F} = \det FF^{- T},$$
we have
\begin{eqnarray}\label{c2}
&&(\det(I + V))_t + u\cdot\nabla(\det(I + V))\\\nonumber
&&= \det(I + V)(I + V)^{- 1}_{ji}\dot{V}_{ij}\\
&&= tr\Bigg\{\dot{V} + \begin{pmatrix}V_{22} & - V_{12}\\\nonumber
    - V_{12} & V_{11}\end{pmatrix}\dot{V}\Bigg\}.
\end{eqnarray}
Return to \eqref{b15}, we in fact have that
\begin{equation}\label{c3}
\dot{V} = D(u) + \frac{1}{2}(\nabla uV + V\nabla u^T) -
\frac{1}{2}(VA - AV)\dot{\theta}.
\end{equation}
By \eqref{c2}, \eqref{c3}, we arrive at
\begin{eqnarray}\label{c4}
&&(\det(I + V))_t + u\cdot\nabla(\det(I + V))\\\nonumber &&=\
trD(u) - \frac{1}{2}tr(VA - AV)\dot{\theta} +
    tr\Bigg\{\nabla uV + \begin{pmatrix}V_{22} & -
    V_{12}\\- V_{12} & V_{11}\end{pmatrix}D(u)\Bigg\}\\\nonumber
&&\quad +\ tr\Bigg\{\begin{pmatrix}V_{22} & - V_{12}\\- V_{12} &
    V_{11}\end{pmatrix}\nabla uV\Bigg\} +
    tr\Bigg\{\begin{pmatrix}V_{22} & - V_{12}\\- V_{12} &
    V_{11}\end{pmatrix}(VA - AV)\Bigg\}.
\end{eqnarray}

We will compute the right side of the above equality term by term.
First of all, with the use of the first equation of \eqref{b}, it
is rather easy to see that
\begin{equation}\label{c5}
trD(u) = 0.
\end{equation}
Making use of \eqref{b10}, we deduce that
\begin{equation}\label{c6}
tr(VA - AV) = tr\begin{pmatrix}2V_{12} & V_{22} - V_{11}\\V_{22} -
V_{11} & - 2V_{12}\end{pmatrix} = 0.
\end{equation}
According to the expression of the Cauchy strain tensor
\eqref{b16}, the third term of \eqref{c4} can be calculated as
\begin{eqnarray}\label{c7}
&&tr\Bigg\{\nabla uV + \begin{pmatrix}V_{22} & - V_{12}\\- V_{12}
& V_{11}\end{pmatrix}D(u)\Bigg\}\\\nonumber &&= \nabla_1u_1V_{11}
+ \nabla_2u_1V_{12} + \nabla_1u_2V_{12} +
\nabla_2u_2V_{22}\\\nonumber &&\quad +\ V_{22}\nabla_1u_1 -
V_{12}(\nabla_1u_2 + \nabla_2u_1) + V_{11}\nabla_2u_2\\\nonumber
&&= trV\nabla\cdot u = 0.
\end{eqnarray}
Similarly, by a straightforward computation, we find that
\begin{eqnarray}\label{c8}
&&tr\Bigg\{\begin{pmatrix}V_{22} & - V_{12}\\- V_{12} &
       V_{11}\end{pmatrix}\nabla uV\Bigg\}\\\nonumber
&&= tr\Bigg\{\begin{pmatrix}V_{22} & - V_{12}\\- V_{12} &
       V_{11}\end{pmatrix}\begin{pmatrix}\nabla_1u_1V_{11} +
       \nabla_2u_1V_{12} & \nabla_1u_1V_{12} + \nabla_2u_1V_{22}\\
       \nabla_1u_2V_{11} + \nabla_2u_2V_{12}
       & \nabla_1u_2V_{12} + \nabla_2u_2V_{22}\end{pmatrix}\Bigg\}\\\nonumber
&&= V_{22}(\nabla_1u_1V_{11} + \nabla_2u_1V_{12}) +
      V_{11}(\nabla_1u_2V_{12} + \nabla_2u_2V_{22})\\\nonumber
&&\quad -\ V_{12}[(\nabla_1u_2V_{11} + \nabla_2u_2V_{12}) +
      (\nabla_1u_1V_{12} + \nabla_2u_1V_{22})]\\\nonumber
&&= V_{11}V_{22}\nabla\cdot u + V_{12}(V_{22}\nabla_2u_1 +
      V_{11}\nabla_1u_2)\\\nonumber
&&\quad -\ V_{12}[\nabla_1u_2V_{11} + \nabla_2u_1V_{22} +
      V_{12}\nabla\cdot u]\\\nonumber &&= 0.
\end{eqnarray}
It remains to deal with the last term of \eqref{c4}. Recall
\eqref{b10} once more, it follows from a simple calculation that
\begin{eqnarray}\label{c9}
&&tr\Bigg\{\begin{pmatrix}V_{22} & - V_{12}\\- V_{12} &
V_{11}\end{pmatrix}(VA - AV)\Bigg\}\\\nonumber &&=
tr\Bigg\{\begin{pmatrix}V_{22} & - V_{12}\\- V_{12} &
V_{11}\end{pmatrix}\begin{pmatrix}2V_{12} & V_{22} - V_{11}\\
V_{22} - V_{11} & - 2V_{12}\end{pmatrix}\Bigg\}\\\nonumber &&=
2V_{22}V_{12} - 2V_{12}(V_{22} - V_{11}) - 2V_{11}V_{12} = 0.
\end{eqnarray}

Finally, combining the above equalities \eqref{c4} with
\eqref{c5}-\eqref{c9}, we arrive at
\begin{equation}\label{c10}
(\det(I + V))_t + u\cdot\nabla(\det(I + V)) = 0.
\end{equation}
This completes the proof Lemma \ref{lem1}.

\begin{rem}\label{rem2}
The incompressibility can be exactly represented as
\begin{equation}\label{c11}
\det F = 1.
\end{equation} The usual incompressible condition $\nabla\cdot u = 0$, the
first equation in \eqref{b}, is the direct consequence of this
identity. On the other hand, $\det(I + V) = 1$ is also a direct
conclusion of \eqref{c11} since $R$ is an orthogonal matrix.
However, the above Lemma illustrates the incompressible
consistency of the the system \eqref{b}.
\end{rem}

To show the global existence of the classical small strain
solutions to system \eqref{b}, we still need to show that the
pressure is a high order term, which will be seen more clear in
the next section. To do so, we will use the following lemma to
find out the weak dissipation of $\nabla\theta$, which will play
an important role for our proof below.

\begin{lem}\label{lem2}
Assume that the third equality of \eqref{bbb} is satisfied. Then
any solution $(u, V)$ of the system \eqref{b}-\eqref{bb} will
satisfy the following identity
\begin{equation}\label{d1}
\nabla\cdot V = A(I+ V)\nabla\theta
\end{equation}
for all the latter time $t \geq 0$.
\end{lem}

$Proof.$ We rewrite the third equation of \eqref{b} as
\begin{equation}\label{d01}
V_t + u\cdot\nabla V = D(u) + \frac{1}{2}(\nabla uV +
      V\nabla u^T) - \frac{1}{2}(VA - AV)\dot{\theta},
\end{equation}
and then apply $\nabla\cdot$ to the resulting equation \eqref{d01}
to yield
\begin{eqnarray}\label{d2}
&&\nabla\cdot V_t + u\cdot\nabla\nabla\cdot V\\\nonumber &&=
\frac{1}{2}\Delta u + \frac{1}{2}AV\nabla\dot{\theta}
    - \begin{pmatrix}\nabla_ju\cdot\nabla V_{1j}\\
    \nabla_ju\cdot\nabla V_{2j}\end{pmatrix}
    +  \frac{1}{2}\nabla\cdot(\nabla uV + V\nabla u^T)\\\nonumber
&&\quad -\ \frac{1}{2}VA\nabla\dot{\theta} -
    \frac{1}{2}\nabla\cdot(VA - AV)\dot{\theta}\\\nonumber
&&= \frac{1}{2}\Delta u + \frac{1}{2}AV\nabla\dot{\theta}
    - \begin{pmatrix}\nabla_ju\cdot\nabla V_{1j}\\
    \nabla_ju\cdot\nabla V_{2j}\end{pmatrix}
    +  \frac{1}{2}\nabla\cdot(\nabla uV + V\nabla
    u^T)\\\nonumber
&&\quad-\ \frac{1}{2}VA\nabla\dot{\theta} -
    \frac{1}{2}\nabla\cdot(VA + AV)\dot{\theta} + A\nabla\cdot
V\dot{\theta}.
\end{eqnarray}
By the first equation of \eqref{b} and a straightforward
computation, we obtain
\begin{equation}\label{d3}
- \begin{pmatrix}\nabla_ju\cdot\nabla V_{1j}\\
    \nabla_ju\cdot\nabla V_{2j}\end{pmatrix} +  \frac{1}{2}\nabla\cdot(\nabla uV + V\nabla u^T)
    = \frac{1}{2}\nabla\cdot(\nabla uV - V\nabla u^T).
\end{equation}
In view of \eqref{b13}, we can rewrite the last term of \eqref{d2}
as
\begin{equation}\label{d02}
A\nabla\cdot V\dot{\theta} = [\omega(u) + \gamma A]\nabla\cdot V.
\end{equation}
Thus, it follows from \eqref{d2}, \eqref{d3} and \eqref{d02} that
\begin{eqnarray}\label{d4}
&&\nabla\cdot V_t + u\cdot\nabla\nabla\cdot V\\\nonumber &&=
\frac{1}{2}\Delta u + \frac{1}{2}AV\nabla\dot{\theta} +
    [\omega(u) + \gamma A]\nabla\cdot V - \frac{1}{2}VA
    \nabla\dot{\theta} \\\nonumber
&&\quad +\ \frac{1}{2}\nabla\cdot(\nabla uV - V\nabla u^T) -
    \frac{1}{2}\nabla\cdot(VA + AV)\dot{\theta}.
\end{eqnarray}

On the other hand, combining \eqref{b21} with \eqref{d01}, and
splitting some terms of the resulting expression, we conclude that
\begin{eqnarray}\label{d03}
 &&\partial_t[A(I + V)\nabla\theta] + u\cdot\nabla[A(I +
    V)\nabla\theta]\\\nonumber
&&= AD(u)\nabla\theta + \frac{1}{2}A(\nabla uV + V\nabla
    u^T)\nabla\theta - \frac{1}{2}(AVA + V)\dot{\theta}\nabla\theta\\\nonumber
    &&\quad +\ A(I + V)\nabla\dot{\theta} - A(I + V)\begin{pmatrix}\nabla_1u\cdot\nabla\theta\\
    \nabla_2u\cdot\nabla\theta\end{pmatrix}\\\nonumber
&&= A\nabla\dot{\theta} + \frac{1}{2}AV\nabla\dot{\theta} +
    \frac{1}{2}AV\nabla\dot{\theta} +
    AD(u)\nabla\theta + \frac{1}{2}A(\nabla uV + V\nabla u^T)\nabla\theta\\\nonumber
&&\quad -\ \frac{1}{2}(AVA + V)\dot{\theta}\nabla\theta - A\begin{pmatrix}\nabla_1u\cdot\nabla\theta\\
    \nabla_2u\cdot\nabla\theta\end{pmatrix} - AV\begin{pmatrix}\nabla_1u\cdot\nabla\theta\\
    \nabla_2u\cdot\nabla\theta\end{pmatrix}.
\end{eqnarray}
In view of \eqref{b13} and the first equation of system \eqref{b},
we obtain, after a straightforward calculation,
\begin{equation}\label{d04}
A\nabla\dot{\theta} = \frac{1}{2}\Delta u + A\nabla\gamma.
\end{equation}
It is rather easy to see that
\begin{equation}\label{d05}
AD(u)\nabla\theta - A\begin{pmatrix}\nabla_1u\cdot\nabla\theta\\
\nabla_2u\cdot\nabla\theta\end{pmatrix} = \omega(u)A\nabla\theta
\end{equation}
and
\begin{equation}\label{d06}
\frac{1}{2}A(\nabla uV + V\nabla u^T)\nabla\theta
- AV\begin{pmatrix}\nabla_1u\cdot\nabla\theta\\
\nabla_2u\cdot\nabla\theta\end{pmatrix} = \frac{1}{2}A(\nabla uV -
V\nabla u^T)\nabla\theta.
\end{equation}
Using \eqref{b11} and \eqref{b13}, we find that
\begin{eqnarray}\label{d07}
&&- \frac{1}{2}(AVA +
   V)\dot{\theta}\nabla\theta\\\nonumber
&&= - \frac{1}{2}A(VA + AV)\dot{\theta}\nabla\theta - (-
\omega_{12}(u) + \gamma)V\nabla\theta\\\nonumber &&= -
\frac{1}{2}A(VA + AV)\dot{\theta}\nabla\theta +
   A\omega(u)V\nabla\theta - \gamma V\nabla\theta.
\end{eqnarray}
Combining \eqref{d03}-\eqref{d07}, we have
\begin{eqnarray}\label{d08}
 &&\partial_t[A(I + V)\nabla\theta] + u\cdot\nabla[A(I +
    V)\nabla\theta]\\\nonumber
&&= \frac{1}{2}\Delta u + \frac{1}{2}AV\nabla\dot{\theta} +
    \omega(u)A(I + V)\nabla\theta - \frac{1}{2}A(VA
    + AV)\dot{\theta}\nabla\theta\\\nonumber
&&\quad +\ A\nabla\gamma + \frac{1}{2}A(\nabla uV - V\nabla
    u^T)\nabla\theta + \frac{1}{2}AV\nabla\dot{\theta} - \gamma
    V\nabla\theta
\end{eqnarray}

Now, we subtract \eqref{d08} from \eqref{d4}, and rearrange the
right side terms of resulting equation. The outcome of this
straightforward calculations is
\begin{eqnarray}\label{d6}
 &&\partial_t[\nabla\cdot V - A(I + V)\nabla\theta] +
u\cdot\nabla[\nabla\cdot V - A(I +
    V)\nabla\theta]\\\nonumber
&&= \omega(u)[\nabla\cdot V - A(I + V)\nabla\theta] -
    \frac{1}{2}(VA + AV)\nabla\dot{\theta}
    - \frac{1}{2}\nabla\cdot(VA + AV)\dot{\theta}\\\nonumber
&&\quad +\ \frac{1}{2}\nabla\cdot(\nabla uV - V\nabla u^T) -
    A\nabla\gamma + \gamma A\nabla\cdot V\\\nonumber
&&\quad +\ [\frac{1}{2}A(VA + AV)\dot{\theta} -
\frac{1}{2}A(\nabla uV -
    V\nabla u^T) + \gamma V]\nabla\theta\\\nonumber
&&= \omega(u)[\nabla\cdot V - A(I + V)\nabla\theta] -
    \frac{1}{2}\nabla\cdot[(VA + AV)\dot{\theta}]\\\nonumber
&&\quad +\ \frac{1}{2}\nabla\cdot(\nabla uV - V\nabla u^T) -
    A\nabla\gamma + \gamma A\nabla\cdot V\\\nonumber
&&\quad +\ [\frac{1}{2}A(VA + AV)\dot{\theta} -
\frac{1}{2}A(\nabla uV -
    V\nabla u^T) + \gamma V]\nabla\theta.
\end{eqnarray}
Our objective is to show that the expression  $\nabla\cdot V - A(I
+ V)\nabla\theta$  satisfies some differential equation whose
solution can be uniquely determined by the initial data. Recall
\eqref{b11} and \eqref{b13}, we get
\begin{eqnarray}
\nonumber && -\ \frac{1}{2}\nabla\cdot[(VA + AV)\dot{\theta}] +
    \frac{1}{2}\nabla\cdot(\nabla uV - V\nabla u^T) -
    A\nabla\gamma + \gamma A\nabla\cdot V\\\nonumber
&&= -\ \frac{1}{2}\nabla\cdot[trVA(- \omega_{12}(u) + \gamma)] -
    A\nabla\gamma + \gamma A\nabla\cdot V\\\nonumber
&& \quad +\ \frac{1}{2}\nabla\cdot(\nabla uV - V\nabla
    u^T)\\\nonumber &&= \frac{1}{2}A\nabla[trV\omega_{12}(u)] -
    \frac{1}{2}A\nabla[(trV + 2)\gamma)] + \gamma A\nabla\cdot V\\\nonumber
&&\quad -\ \frac{1}{2}\nabla\cdot[A(\nabla_ku_1V_{k2}
    \quad -\ \nabla_ku_2V_{k1})]\\\nonumber
&&= \frac{1}{2}A\nabla[trV\omega_{12}(u) - (\nabla_ku_1V_{k2} -
    \nabla_ku_2V_{k1})] + \gamma A\nabla\cdot V\\\nonumber
&&\quad -\ \frac{1}{2}A\nabla[(trV + 2)\gamma)].
\end{eqnarray}
Thus, by \eqref{b02}, we can easily deduce that
\begin{eqnarray}\label{d7}
&&- \frac{1}{2}\nabla\cdot[(VA + AV)\dot{\theta}] +
\frac{1}{2}\nabla\cdot(\nabla uV - V\nabla u^T)\\\nonumber && -\
A\nabla\gamma + \gamma A\nabla\cdot V = \gamma A\nabla\cdot V.
\end{eqnarray}
On the other hand, it follows from \eqref{b11}, \eqref{b02} and
\eqref{b13} that
\begin{eqnarray}\label{d8}
 &&\frac{1}{2}A(VA + AV)\dot{\theta} -
    \frac{1}{2}A(\nabla uV - V\nabla u^T) + \gamma V\\\nonumber
&&= - \gamma AA(I +V) - \gamma I - \frac{1}{2}trV[-
    \omega_{12}(u)\\\nonumber
&&\quad +\ \gamma]I - \frac{1}{2}A(\nabla uV - V\nabla
    u^T)\\\nonumber
&&= - \gamma AA(I +V) - \frac{1}{2}[- trV\omega_{12}(u) +
    \gamma(trV + 2)]I\\\nonumber
&&\quad -\ \frac{1}{2}A(\nabla uV - V\nabla u^T)\\\nonumber &&= -
\gamma AA(I +V) - \frac{1}{2}(\nabla_ku_1V_{k2} -
    \nabla_ku_2V_{k1})I\\\nonumber
&&\quad +\ \frac{1}{2}A(\nabla uV + V\nabla u^T)\\\nonumber &&= -
\gamma AA(I +V).
\end{eqnarray}
Hence we conclude from \eqref{d6}-\eqref{d8} that
\begin{eqnarray}\label{d9}
&&\partial_t[\nabla\cdot V - A(I + V)\nabla\theta] +
u\cdot\nabla[\nabla\cdot V - A(I +
    V)\nabla\theta]\\\nonumber
&&= [\omega(u) + \gamma A][\nabla\cdot V - A(I + V)\nabla\theta].
\end{eqnarray}
which concludes the proof of the lemma, since the above quantity
will remain zero all the time with zero initial data.

\begin{rem}\label{rem3}
In order to carry out the mechanical background of the above
lemma, we again go back to the deformation tensor $F$. By its
definition, the following identity
\begin{equation}\nonumber
\nabla\cdot F^T = 0
\end{equation}
is proved in \cite{Lei4, Lei5}. Insert \eqref{b1} into the above
equality, we find
$$\partial_iR_{ij} + R_{kj}\partial_iV_{ik} + V_{ik}\partial_iR_{kj} = 0$$
which implies
\begin{equation}\nonumber
\partial_iV_{il} + R_{lj}\partial_iR_{ij} + R_{lj}V_{ik}\partial_iR_{kj} = 0
\end{equation}
since $R$ is an orthogonal matrix. But then employing \eqref{b3} a
straightforward calculation yields \eqref{d1}.
\end{rem}

\section{Proof of Theorem \ref{thm1}}

In this section we prove Theorem \ref{thm1}. Weak dissipations on
the strain $V$ and the gradient of rotation angle $\nabla\theta$
are found by introducing auxiliary functions $w$ and $\Theta$
below. The way of defining such functions reveals the intrinsic
dissipative nature of the system. In what follows, $\|\cdot\|$ is
used to denote the standard $L^2(R^2)$ norm, $<\cdot, \cdot>$ and
$(\cdot, \cdot)$ the $R^d$ and $L^2(R^2)^d$ inner products with $d
= 1,\ 2$ or 4. The proof is divided into four steps.

\bigskip
\textbf{Step 1. Standard Energy Estimates.}
\bigskip

By taking the $L^2$ inner product of the second and third
equations of system \eqref{b} with $u$ and $V$, and then adding up
the resulting equations together, we obtain
\begin{eqnarray}
\nonumber&&\frac{1}{2}\frac{d}{dt}\int_{R^2}(|u|^2 + |V|^2)\ dx +
  \mu\|\nabla u\|^2\\\nonumber
&& =- \int_{R^2}u\cdot\nabla\frac{|u|^2 + |V|^2}{2} + <u, \nabla
  p>\ dx\\\nonumber
&&\quad +\ \int_{R^2}<V, \frac{1}{2}(\nabla uV + V\nabla u^T)
  + \frac{1}{2}\omega_{12}(u)(VA - AV)>\ dx\\\nonumber
&& \quad +\ \int_{R^2}<u, 2\nabla\cdot V > + <V, D(u)>\
  dx\\\nonumber
&&\quad +\ \int_{R^2}<u, \nabla\cdot(VV^T)>\ dx - \int_{R^2}<V,
  \frac{1}{2}\gamma(VA - AV)>\ dx.
\end{eqnarray}
On the other hand, the fourth equation of system \eqref{b} gives
\begin{eqnarray}
\nonumber&&\frac{1}{2}\frac{d}{dt}\int_{R^2}2trV\ dx  = -
  \int_{R^2}u\cdot\nabla trV\ dx + \int_{R^2}trD(u)\ dx\\\nonumber
&&\quad +\ \int_{R^2}tr\big(\nabla uV + \frac{1}{2}(VA-
  AV)(\omega_{12}(u) - \gamma)\big)\ dx.
\end{eqnarray}
Adding up the above two equation together, we have
\begin{eqnarray}\label{f1}
&&\frac{1}{2}\frac{d}{dt}\int_{R^2}(|u|^2 + |V|^2 +
  2trV)\ dx + \mu\|\nabla u\|^2\\
\nonumber&& = - \int_{R^2}u\cdot\nabla\big(\frac{|u|^2 + |V|^2}{2}
  + trV\big) + <u, \nabla p>\ dx + \int_{R^2}trD(u)\ dx\\\nonumber
&&\quad +\ \int_{R^2}<V, \frac{1}{2}(\nabla uV + V\nabla u^T)
  + \frac{1}{2}\omega_{12}(u)(VA - AV)>\ dx\\\nonumber
&& \quad +\ \int_{R^2}<u, 2\nabla\cdot V > + <V, D(u)>\
  dx + \int_{R^2}tr(\nabla uV)\ dx \nonumber\\\nonumber
&&\quad +\ \int_{R^2}<u, \nabla\cdot(VV^T)>\ dx -
  \int_{R^2}<V, \frac{1}{2}\gamma(VA - AV)>\ dx\\\nonumber
&&\quad+\ \int_{R^2}tr\frac{1}{2}\big((VA- AV) (\omega_{12}(u) -
  \gamma)\big)\ dx.
\end{eqnarray}

Noting that $V$ is symmetric, and with the aid of the first
equation of \eqref{b}, it is rather easy to see that
\begin{equation}\label{f2}
\begin{cases}
- \int_{R^2}u\cdot\frac{1}{2}(|u|^2 + |V|^2 + 2trV) + <u, \nabla
  p> dx + \int_{R^2}trD(u)\ dx = 0,\\
\int_{R^2}<u, 2\nabla\cdot V > + <V, D(u)>\ dx +
  \int_{R^2}tr(\nabla uV)\ dx = 0.
\end{cases}
\end{equation}
Noting the definition in \eqref{b10} and \eqref{b02}, one has
\begin{eqnarray}\label{f3}
&&\int_{R^2}<V, \frac{1}{2}(\nabla uV + V\nabla
  u^T) + \frac{1}{2}\omega_{12}(u)(VA - AV)>\ dx\\\nonumber
&&\quad +\ \int_{R^2}<u, \nabla\cdot(VV^T)>\ dx - \int_{R^2}<V,
  \frac{1}{2}\gamma(VA - AV)>\ dx\\\nonumber
&&\quad +\ \int_{R^2}tr\frac{1}{2}\big((VA- AV) (\omega_{12}(u) -
  \gamma)\big)\ dx\\\nonumber
&&= \int_{R^2}<V, \nabla uV> + <u, \nabla\cdot(VV^T)>\
  dx\\\nonumber
&&\quad +\ \int_{R^2}\frac{1}{2}(\omega_{12}(u)
  - \gamma)\big(tr(VA- AV) + <V, VA - AV>\big)\ dx\\\nonumber
&&= 0.
\end{eqnarray}
Then, combining \eqref{f1}, \eqref{f2} with \eqref{f3}, we have
\begin{equation}\nonumber
\frac{1}{2}\frac{d}{dt}\int_{R^2}(|u|^2 + |V|^2 + 2trV)\ dx +
\mu\|\nabla u\|^2 = 0,
\end{equation}
Finally, by using \eqref{coreq}, we arrive at
\begin{equation}\label{f4}
\frac{1}{2}\frac{d}{dt}\int_{R^2}(|u|^2 + |V_{11} - V_{22}|^2 +
|V_{12} + V_{21}|^2)\ dx + \mu\|\nabla u\|^2 = 0.
\end{equation}

By taking the $L^2$ inner product of the second and third
equations of system \eqref{b} with $u$ and $2V$, and then adding
up the resulting equations together, we obtain
\begin{eqnarray}
\nonumber&&\frac{1}{2}\frac{d}{dt}\int_{R^2}|u|^2 +
    2|V|^2dx + \int_{R^2}u\cdot\nabla(\frac{|u|^2}{2}
    + |V|^2)dx + (u, \nabla p)\\\nonumber
&&= \Big(V, \nabla uV + V\nabla u^T +
    \omega_{12}(u)(VA - AV) - \gamma(VA - AV)\Big)\\\nonumber
&&\ \ \ \ +\ \big(u, \nabla\cdot(VV^T)\big) + (u, 2\nabla\cdot V)
+ \big(2V, D(u)\big) +
    \mu(u, \Delta u),
\end{eqnarray}
which gives that
\begin{eqnarray}\label{f5}
&&\frac{1}{2}\frac{d}{dt}\int_{R^2}|u|^2 + 2|V|^2dx +
  \mu\|\nabla u\|^2\\\nonumber
&&\leq\ \|\nabla u\|_{L^\infty}\|V\|^2\\\nonumber &&\leq
  C\|\nabla V\|\|\Delta V\|(\|\Delta\nabla u\| + \|\nabla u\|)\\\nonumber
&&\leq C\mu\|\nabla V\|_{L^2}\big(\|\nabla u\|^2 + \|\Delta\nabla
u\|^2 +
  \frac{1}{\mu^2}\|\Delta V\|\big).
\end{eqnarray}

In order to exploit the higher order energy estimates, we apply
$\Delta$ to the second and third equations of system \eqref{b},
and then take $L^2$ inner product of the resulting equations with
$\Delta u$ and $2\Delta V$, respectively, to yield
\begin{eqnarray}\label{f6}
&&\frac{1}{2}\frac{d}{dt}\int_{R^2}|\Delta u|^2 +
    2|\Delta V|^2dx + (\Delta u, \nabla\Delta p)\\\nonumber
&&\quad +\ \int_{R^2}u\cdot\nabla(\frac{|\Delta u|^2}{2}
    + |\Delta V|^2)dx\\\nonumber
&&=\ \Big(\Delta V, \Delta\big[\nabla uV + V\nabla u^T +
    \omega_{12}(u)(VA - AV)\\\nonumber
&&\quad -\ \gamma(VA - AV)\big]\Big) + \big(\Delta u,
    \Delta\nabla\cdot(VV^T)\big) + (\Delta u, 2\Delta\nabla\cdot V)\\\nonumber
&&\quad +\ \big(2\Delta V, \Delta D(u)\big) +
    \mu(\Delta u, \Delta\Delta u) - \Big(\Delta u, \big[\Delta(u\cdot\nabla u) -
    u\cdot\nabla\Delta u\big]\Big)\\\nonumber
&&\quad -\ \Big(2\Delta V, \big[\Delta(u\cdot\nabla V) -
    u\cdot\nabla\Delta V\big]\Big).
\end{eqnarray}
Applying similar arguments as before, we conclude that
\begin{equation}\label{f7}
\begin{cases}
\int_{R^2}u\cdot\nabla(\frac{|\Delta u|^2}{2} + |\Delta V|^2)dx +
(\Delta u, \nabla\Delta p)\\
\quad -\ \mu(\Delta u, \Delta\Delta u) =
\mu\|\Delta\nabla u\|^2,\\
(\Delta u, 2\Delta\nabla\cdot V) + \big(2\Delta V, \Delta
D(u)\big) = 0.
\end{cases}
\end{equation}

To estimate the rest terms of the right side of \eqref{f6}, and
also for later use, we need the following Lemmas. Their proofs can
found in many literatures, see \cite{Lei4}, for example.

\begin{lem}\label{lem3}
Assume $v \in W^{k,2}(R^2), k \geq 3$. The following interpolation
inequalities hold.
\begin{enumerate}
     \item For $1 \leq s \leq k$,
          \begin{eqnarray}\nonumber
               \|v\|_{L^{4}} &\leq& C\|v\|^{1 - \frac{1}{2s}}
                   \|\nabla^{s} v\|^{\frac{1}{2s}},\\\nonumber
               \|\nabla v\|_{L^{4}} &\leq& C\|v\|^{1 - \frac{3}{2(s +
                   1)}} \|\nabla^{s}\nabla v\|^{\frac{3}{2(s + 1)}},\\\nonumber
               \|\Delta v\|_{L^{4}} &\leq& C\|v\|^{1 - \frac{5}{2(s +
2)}}\|\nabla^s\Delta v\|^{\frac{5}{2(s + 2)}}.
         \end{eqnarray}
     \item For $2 \leq s \leq k$,
          \begin{eqnarray}\nonumber
               \|v\|_{L^{\infty}} &\leq& C\|v\|^{1 - \frac{1}{s}}
                   \|\nabla^{s} v\|^{\frac{1}{s}},\\\nonumber
               \|\nabla v\|_{L^{\infty}} &\leq& C\|v\|^{1 - \frac{2}{s
                   + 1}}\|\nabla^s\nabla v\|^{\frac{2}{s + 1}}.
          \end{eqnarray}
\end{enumerate}
\end{lem}

The following Lemma will also be used for several times below.

\begin{lem}\label{lem4}
Assume $f,\ g \in H^{s}(R^2)$ for some $s > 0$. Then there exists
a universal constant $C$ depending only on $s$ such that
$$\|\nabla^s(fg)\| \leq C(\|f\|_{L^\infty}\|\nabla^s g\| +
\|g\|_{L^\infty}\|\nabla^s f\|).$$
\end{lem}

Now we estimate the second line of \eqref{f6}. Invoking the
definition of $\gamma$ in \eqref{b02},  and with the aid of Lemma
\ref{lem4}, we deduce that
\begin{eqnarray}\label{f8}
&&\Big|\Big(\Delta V, \Delta\big[\nabla uV + V\nabla u^T
    + \omega_{12}(u)(VA - AV) - \gamma(VA - AV)\big]\Big)\Big|\\\nonumber
&&\leq\ C\|\Delta V\|\Big(\|\Delta V\|\|\nabla u\|_{L^\infty} +
    \|\nabla\Delta u\|\|V\|_{L^\infty} + \|\nabla\Delta
    u\|\|V\|_{L^\infty}^2\\\nonumber
&&\ \ \ \ +\ \|\nabla u\|_{L^\infty}\big[\|V\|_{L^\infty}\|\Delta
    V\| + \|V\|_{L^\infty}^2\|\Delta\frac{1}{2 + trV}\|\big]\Big)\\\nonumber
&&\leq\ C\mu\|V\|_{H^2}(\|\nabla u\|^2 + \|\nabla\Delta u\|^2 +
    \frac{1}{\mu^2}\|\Delta V\|^2)
\end{eqnarray}
provided that
\begin{equation}\label{z1}
\|V\|_{L^\infty} < 1,
\end{equation}
which will be verified below.

Similarly, it is easy to show that
\begin{eqnarray}\label{f9}
&&\big|\big(\Delta u, \Delta\nabla\cdot(VV^T)\big)\big|\\\nonumber
&&\leq \|V\|_{L^\infty}\|\Delta V\|\|\nabla\Delta u\|\\\nonumber
&&\leq C\mu\|V\|_{H^2}(\|\nabla\Delta u\|^2 +
  \frac{1}{\mu^2}\|\Delta V\|^2).
\end{eqnarray}

At last, it remains to estimate the last line of \eqref{f6}.
Again, by Sobolev imbedding theorem and Lemma \ref{lem4}, it leads
to
\begin{eqnarray}\label{f10}
&&\Big|- \Big(\Delta u, \big[\Delta(u\cdot\nabla u) -
    u\cdot\nabla\Delta u\big]\Big) - \Big(2\Delta V, \big[\Delta(u\cdot\nabla V)
    - u\cdot\nabla\Delta V\big]\Big)\Big|\\\nonumber
&&\leq\ C\|\Delta u\|\big(\|\nabla u\|_{L^\infty}\|\Delta u\| +
    \|\nabla\nabla u\|\|u\|_{L^\infty}\big)\\\nonumber
&&\ \ \ \ +\ C\|\Delta V\|\big(\|\nabla u\|_{L^\infty}\|\Delta V\|
    + \|\nabla\nabla u\|\|V\|_{L^\infty}\big)\\\nonumber
&&\leq\ C(1 + \mu)(\|u\|_{H^2} + \|V\|_{H^2})(\|\nabla u\|^2 +
    \|\nabla\Delta u\|^2 + \frac{1}{\mu^2}\|\Delta V\|^2).
\end{eqnarray}

Combining all these estimates \eqref{f7}-\eqref{f10} with
\eqref{f5}, and then noting \eqref{f6}, we can finally conclude
that
\begin{eqnarray}\label{f}
&&\frac{1}{2}\frac{d}{dt}\int_{R^2}\big(|u|^2 + |\Delta
    u|^2 + 2|V|^2 + 2|\Delta V|^2\big)dx + \mu\big(\|\nabla u\|^2
    + \|\nabla\Delta u\|^2\big)\\\nonumber
&&\leq\ C(1 + \mu)(\|u\|_{H^2} + \|V\|_{H^2})\big(\|\nabla u\|^2 +
    \|\nabla\Delta u\|^2 + \frac{1}{\mu^2}\|\Delta V\|^2\big)
\end{eqnarray}
provided \eqref{z1} holds.

\bigskip
\textbf{Step 2. Weak Dissipation for The Strain Matrix
$V$.}
\bigskip

To get dissipative energy estimates, by \eqref{f}, it is clear
that we should explore the dissipation of $V$.

Let
\begin{equation}\label{g1}
w = \Delta u + \frac{2}{\mu}\nabla\cdot V.
\end{equation}

Now we apply $\Delta$ to the second equation of system \eqref{b},
and respectively, $\frac{1}{\mu}\nabla\cdot$ to the third one, and
then adding up the outcomes to get a new equation. In terms of $w$
introduced in \eqref{g1}, this new equation can be written as
\begin{eqnarray}\label{g2}
w_t + u\cdot\nabla w + \nabla\Delta p = \mu\Delta w +
\frac{1}{\mu}\Delta u + f.
\end{eqnarray}
where
\begin{eqnarray}\label{g3}
 &&f = - \big[\Delta(u\cdot\nabla u) - u\cdot\nabla\Delta
      u\big] - \frac{2}{\mu}\big[\nabla\cdot(u\cdot\nabla V)\\\nonumber
 &&\quad -\ u\cdot\nabla(\nabla\cdot V)\big] + \Delta\nabla\cdot(VV^T) +
      \frac{1}{\mu}\nabla\cdot\big[(\nabla uV + V\nabla u^T)\\\nonumber
 &&\quad +\ \omega_{12}(u)(VA - AV) - \gamma(VA - AV)\big].
\end{eqnarray}

Take the inner product of \eqref{g2} with $w$ in $(L^2(R^2))^2$,
and use the similar arguments as in \eqref{f2}, we have
\begin{equation}\label{g4}
\frac{1}{2}\frac{d}{dt}\int_{R^2}|w|^2dx + \mu\|\nabla w\|^2 =
(\Delta p, \nabla\cdot w) + \frac{1}{\mu}(\Delta u, w) + (f, w).
\end{equation}

As before, we should estimate the right side of \eqref{g4} term by
term. Firstly, it is obvious that
\begin{equation}\label{g5}
\frac{1}{\mu}|(\Delta u, w)| \leq \frac{\mu}{4}\|\nabla w\|^2 +
\frac{1}{\mu^3}\|\nabla u\|^2.
\end{equation}

On the other hand, noting that $\nabla\cdot u = 0$, we can use
integration by parts and Lemma \ref{lem4} to estimate the last
term of the right side of \eqref{g4} as follows:
\begin{eqnarray}\label{g6}
&&|(f, w)| \leq \Big|\Big(\big[\Delta(u_iu_j) -
    u_i\Delta u_j\big], \nabla_i w_j\Big)\Big| + \big|\big(\Delta(VV^T),
    \nabla w\big)\big|\\\nonumber
&&\quad +\ \frac{1}{\mu}\Big|\Big(\big[(\nabla uV + V\nabla u^T) +
    \omega_{12}(u)(VA - AV) - \gamma(VA - AV)\big],\\\nonumber
&&\quad \nabla w\Big)\Big| +
    \frac{2}{\mu}\Big|\Big(\big[\nabla_k\cdot(u_i V_{jk}) -
    u_i(\nabla_kV_{jk})\big], \nabla_iw_j\Big)\Big|\\\nonumber
&&\leq C\|\nabla w\|\Big[\big(\|u\|_{L^\infty}\|\Delta u\| +
    \|\nabla u\|_{L^\infty}\|\nabla u\|\big)\\\nonumber
&&\quad +\ \frac{1}{\mu}\|V\|_{L^\infty}\|\nabla u\|+\
  \|V\|_{L^\infty}\|\Delta V\|
  + \frac{1}{\mu}\|V\|_{L^\infty}^2\|\nabla u\|\Big]\\\nonumber
&&\leq C(\mu + \frac{1}{\mu})(\|u\|_{H^2} +
    \|V\|_{H^2})\\\nonumber
&&\quad \times\big(\|\nabla u\|^2 + \|\Delta\nabla u\|^2 +
    \|\nabla w\|^2 + \frac{1}{\mu^2}\|\Delta V\|^2\big)
\end{eqnarray}
provided $\|V\|_{L^\infty} < 1$.

It also remains to estimate the first term of the right side of
\eqref{g4}. For this purpose, we apply the divergence operator for
the momentum equation of system \eqref{b} to get
\begin{equation}\label{g7}
\Delta p = - tr(\nabla u\nabla u) + \nabla\cdot[\nabla\cdot
(VV^T)] + 2\nabla\cdot(\nabla\cdot V).
\end{equation}
Thus, by \eqref{g7}, Lemma \ref{lem3}, Lemma \ref{lem4} and
Sobolev imbedding theorem, we can easily deduce that
\begin{eqnarray}\label{g8}
&&|(\Delta p, \nabla w)| \leq  \big|\big(tr(\nabla
    u\nabla u), \nabla w\big)\big|\\\nonumber
&&\quad +\ \Big|\Big(\nabla\cdot[\nabla\cdot
    (VV^T)], \nabla w\Big)\Big| + 2\big|\big(\nabla\cdot(\nabla\cdot V), \nabla
    w\big)\big|\\\nonumber
&&\leq C\|\nabla w\|\Big(\|\nabla u\|_{L^4}^2 +
    \|V\|_{L^\infty}\|\Delta V\|\Big)\\\nonumber
&&\quad +\ \frac{\mu}{4}\|\nabla w\|^2 +
    \frac{1}{\mu}\|\nabla\cdot(\nabla\cdot V)\|^2\\\nonumber
&&\leq C(1 + \mu)(\|u\|_{H^2} + \|V\|_{H^2})\big(\|\nabla u\|^2 +
    \|\nabla\Delta u\|^2\\\nonumber
&&\quad +\ \|\nabla w\|^2 + \frac{1}{\mu^2}\|\Delta V\|^2\big)+
    \frac{\mu}{4}\|\nabla w\|^2 + \frac{1}{\mu}\|\nabla\cdot(\nabla\cdot V)\|^2.
\end{eqnarray}
We must deal with the last term of the above estimates. Now Lemma
\ref{lem2} enters the argument in this stage. Noting the
definition of $\nabla^\perp$, it is easy to see that
$\nabla\cdot\nabla^\perp = 0$. Thus, by Lemma \ref{lem3} and a
straightforward calculation, we have
\begin{eqnarray}\label{g9}
&&\frac{2}{\mu}\|\nabla\cdot(\nabla\cdot V)\|^2 =
\frac{2}{\mu}\|\nabla\cdot(AV\nabla\theta)\|^2\\\nonumber &&\leq
\frac{C}{\mu}\big(\|\nabla V\|_{L^4}\|\nabla\theta\|_{L^4}\big)^2
+ \big(\|V\|_{L^\infty}\|\nabla^2\theta\|\big)^2\\\nonumber &&\leq
\frac{C}{\mu}\|\nabla V\|\|\nabla\theta\|\|\Delta V\|
\|\Delta\theta\| +
\frac{C}{\mu}\|V\|_{L^\infty}^2\|\Delta\theta\|^2\\\nonumber
&&\leq C(1 + \mu)(\|V\|_{H^2} + \|u\|_{H^2})(\|\nabla\Delta u\|^2
+ \frac{1}{\mu^2}\|\Delta\theta\|^2 + \frac{1}{\mu^2}\|\Delta
V\|^2)
\end{eqnarray}
provided \eqref{z1} is satisfied and
\begin{equation}\label{z2}
\|\nabla\theta\| \leq 1
\end{equation}
holds.

Finally, the combination of \eqref{g4}-\eqref{g9} leads to
\begin{eqnarray}\label{g}
&&\frac{d}{dt}\int_{R^2}|w|^2dx + \mu\|\nabla
    w\|^2\\\nonumber
&&\leq\ C(\mu + \frac{1}{\mu})(\|u\|_{H^2} +
    \|V\|_{H^2})\big(\|\nabla u\|^2 + \|\Delta\nabla u\|^2\\\nonumber
&&\quad +\ \|\nabla
    w\|^2 + \frac{1}{\mu^2}\|\Delta V\|^2 + \frac{1}{\mu^2}\|\Delta\theta\|^2\big)
+\ \frac{2}{\mu^3}\|\nabla u\|^2
\end{eqnarray}
provided \eqref{z2} holds.

\bigskip
\textbf{Step 3. Weak Dissipation for $\nabla\theta$.}
\bigskip

At this stage, to extract the dissipative nature of the system, it
is clear that we need also explore some kind of weak dissipation
of $\Delta\theta$. For this purpose, we introduce the auxiliary
variable $\Theta$:
\begin{equation}\label{h1}
\Theta = \Delta u + \frac{2}{\mu}\nabla^\perp\theta.
\end{equation}

Rewrite the momentum equation as follows:
\begin{equation}\label{h2}
u_t + u\cdot\nabla u + \nabla p = \mu\Delta u +
2\nabla^\perp\theta + \nabla\cdot (VV^T)
      + 2AV\nabla\theta.
\end{equation}
Then, apply $\Delta$ to \eqref{h2} and $\frac{2}{\mu}\nabla\cdot$
to the last equation of \eqref{b}, and then add up the resulting
equations together, we obtain a new equation in terms of $\Theta$:
\begin{equation}\label{h3}
\Theta_t + u\cdot\nabla\Theta + \nabla\Delta p = \Delta\Theta +
\frac{1}{\mu}\Delta u + g.
\end{equation}
where
\begin{eqnarray}\label{h4}
 &&g = \Delta\nabla\cdot(VV^T) + 2\Delta(AV\nabla\theta)
     + \frac{2}{\mu}\nabla^\perp\gamma\\
&&\quad -\ \big[\Delta(u\cdot\nabla u) - u\cdot\nabla\Delta u\big]
- \frac{2}{\mu}\begin{pmatrix}-
\nabla_2u\cdot\nabla\theta\\\nonumber
    \nabla_1u\cdot\nabla\theta\end{pmatrix}.
\end{eqnarray}

Similarly, we take the $L^2$ inner product of \eqref{h3} with
$\Theta$ to get
\begin{equation}\label{h5}
\frac{1}{2}\frac{d}{dt}\int_{R^2}|\Theta|^2dx + \mu\|\nabla
\Theta\|^2 = (\Delta p, \nabla\cdot \Theta) + \frac{1}{\mu}(\Delta
u, \Theta) + (g, \Theta).
\end{equation}
By the expression of $\Theta$ in \eqref{h1} and the first equation
of system \eqref{b}, it is easy to see that
\begin{equation}\label{h6}
(\Delta p, \nabla\cdot \Theta) = 0.
\end{equation}
A straightforward computation implies
\begin{equation}\label{h7}
\frac{1}{\mu}\big|(\Delta u, \Theta)\big| \leq
\frac{\mu}{4}\|\nabla\Theta\|^2 + \frac{1}{\mu^3}\|\nabla u\|^2.
\end{equation}

To estimate the last term in \eqref{h5}, we insert \eqref{h4} into
$(g, \Theta)$ and then use integration by parts to yield
\begin{eqnarray}\label{h8}
 &&|(g, \Theta)| \leq \Big|\Big(\big[\Delta(u_iu_j) -
    u_i\Delta u_j\big], \nabla_i \Theta_j\Big)\Big| +
    \frac{2}{\mu}\Big|\Big(\begin{pmatrix}- \nabla_2\theta\nabla_2u\\
    \nabla_1\theta\nabla_1u\end{pmatrix}, \Theta\Big)\Big|\\\nonumber
&&\quad +\ \big|\big(\Delta(VV^T), \nabla\Theta\big)\big| +
    \frac{2}{\mu}\big|\big(\gamma, \nabla^\perp\cdot\Theta\big)\big|
    + 2\big|\big(\Delta(AV\nabla\theta), \Theta\big)\big|.
\end{eqnarray}
By a similar argument as before, we have
\begin{eqnarray}\label{h9}
&&\Big|\Big(\big[\Delta(u_iu_j) -
    u_i\Delta u_j\big], \nabla_i \Theta_j\Big)\Big|
    + \big|\big(\Delta(VV^T), \nabla\Theta\big)\big| +
    2\big|\big(\Delta(AV\nabla\theta), \Theta\big)\big|\\\nonumber
&&\leq\  C\|\nabla\Theta\|\Big[\big(\|u\|_{L^\infty}\|\Delta u\| +
    \|\nabla u\|_{L^\infty}\|\nabla u\|\big) + \|V\|_{L^\infty}\|\Delta V\|\Big]\\\nonumber
&&\quad +\
    C\|\nabla\Theta\|\big(\|V\|_{L^\infty}\|\Delta\theta\|
    + \|\nabla V\|_{L^4}\|\nabla\theta\|_{L^4}\big)\\\nonumber
&&\leq\ C(1 + \mu^{\frac{3}{2}})(\|u\|_{H^2} +
    \|V\|_{H^2} + \frac{1}{\mu}\|\nabla\theta\|)\\\nonumber
&&\quad \times\ \big(\|\nabla u\|^2 + \|\Delta\nabla u\|^2 +
    \|\nabla\Theta\|^2 + \frac{1}{\mu^2}\|\Delta V\|^2 + \frac{1}{\mu^2}\|\Delta\theta\|^2\big),
\end{eqnarray}
where we used Lemma \ref{lem3} and Lemma \ref{lem4}. On the other
hand, by \eqref{b02}, it follows that
\begin{eqnarray}\label{h10}
&&\frac{2}{\mu}\Big|\Big(\begin{pmatrix}- \nabla_2\theta\nabla_2u\\
    \nabla_1\theta\nabla_1u\end{pmatrix}, \Theta\Big)\Big| +
    \frac{2}{\mu}\big|\big(\gamma, \nabla^\perp\cdot\Theta\big)\big|\\\nonumber
&&\leq\ \frac{C}{\mu}\big(\|\nabla
    u\|\|\nabla\theta\|_{L^4}\|\Theta\|_{L^4} +
    \|V\|_{L^\infty}\|\nabla u\|\nabla\Theta\|\big)\\\nonumber
&&\leq\ \frac{C}{\mu}\big(\|\nabla
u\|\|\nabla\theta\|^{\frac{1}{2}}\|\Theta\|^{\frac{1}{2}}\|\Delta
    \theta\|^{\frac{1}{2}}\|\nabla\Theta\|^{\frac{1}{2}}
    + \|V\|_{L^\infty}\|\nabla u\|\nabla\Theta\|\big)\\\nonumber
&&\leq\ C(1 + \frac{1}{\mu})(\|\Theta\| +
    \|V\|_{H^2} + \frac{1}{\mu}\|\nabla\theta\|)\\\nonumber
&&\quad \times\ \big(\|\nabla u\|^2 + \|\nabla\Theta\|^2
    + \frac{1}{\mu^2}\|\Delta\theta\|^2\big).
\end{eqnarray}
The combination of \eqref{h8}-\eqref{h10} implies that
\begin{eqnarray}\label{h11}
 &&|(g, \Theta)| \leq C(\mu^2 +
    \frac{1}{\mu})(\|u\|_{H^2} + \|V\|_{H^2} + \|\Theta\|
    + \frac{1}{\mu}\|\nabla\theta\|)\\\nonumber
&&\quad \times\ \big(\|\nabla u\|^2 + \|\Delta\nabla u\|^2 +
    \|\nabla\Theta\|^2 + \frac{1}{\mu^2}\|\Delta\theta\|^2
    + \frac{1}{\mu^2}\|\Delta V\|^2\big).
\end{eqnarray}

Finally, combining \eqref{h5}-\eqref{h7} with \eqref{h11}, and
noting the definition of $\Theta$ in \eqref{h1}, we arrive at
\begin{eqnarray}\label{h}
\nonumber&&\frac{d}{dt}\int_{R^2}|\Theta|^2dx + \mu\|\nabla
    \Theta\|^2\\\nonumber
&&\leq\ \frac{C}{\mu^3}\|\nabla
    u\|^2 + C(\mu^2 + \frac{1}{\mu})(\|u\|_{H^2} +
    \|V\|_{H^2} + \|\Theta\|)\\\nonumber
&&\quad \times\ \big(\|\nabla u\|^2 + \|\Delta\nabla u\|^2 +
    \|\nabla\Theta\|^2 + \frac{1}{\mu^2}\|\Delta V\|^2\big).
\end{eqnarray}

\textbf{Step 4. A $Priori$ Weakly Dissipative Energy Estimates.}

This step ties everything together. Firstly, by Hodge's
decomposition, we calculate
\begin{eqnarray}\label{i1}
&&\Delta V = \nabla\nabla\cdot V -
    \nabla\times\nabla\times V\\\nonumber
&&= \nabla\nabla\cdot V -
    \begin{pmatrix}\nabla_2(\nabla_1V_{12} - \nabla_2V_{11})
    & - \nabla_1(\nabla_1V_{12} - \nabla_2V_{11})\\
    \nabla_2(\nabla_1V_{22} - \nabla_2V_{21})
    & - \nabla_1(\nabla_1V_{22} - \nabla_2V_{21})\end{pmatrix}\\\nonumber
&&= \nabla\nabla\cdot V + (\nabla^\perp)^2trV -
    \nabla^\perp(A\nabla\cdot V).
\end{eqnarray}
Thus, noting \eqref{coreq}, and by \eqref{g1} and Lemma
\ref{lem4}, we obtain
\begin{eqnarray}\label{i2}
&&\|\Delta V\| \leq C\mu(\|\Delta\nabla u\| + \|\nabla
    w\| + \|\Delta trV\|)\\\nonumber
&&\leq C\mu(\|\Delta\nabla u\| + \|\nabla w\|) + C\mu\|\Delta
    V\|\|V\|_{L^\infty}
\end{eqnarray}
which implies that
\begin{equation}\label{i3}
\|\Delta V\| \leq C\mu(\|\Delta\nabla u\| + \|\nabla w\|)
\end{equation}
provided that
\begin{equation}\label{z3}
\|V\|_{L^\infty} \leq \frac{1}{2\mu C}.
\end{equation}
On the other hand, by \eqref{h1},
\begin{equation}\label{i4}
\begin{cases}
\|\nabla\theta\| \leq C\mu(\|\Delta u\| + \|\Theta\|),\ \ \ \ \ \\
\|\Delta\theta\| \leq C\mu(\|\Delta\nabla u\| + \|\nabla\Theta\|).
\end{cases}
\end{equation}

Now we can proceed with the dissipative energy method. Plug
\eqref{i4} into \eqref{f}, we find that
\begin{eqnarray}\label{ff}
&&\frac{1}{2}\frac{d}{dt}\int_{R^2}\big(|u|^2 + |\Delta
    u|^2 + 2|V|^2 + 2|\Delta V|^2\big)dx
    + \mu\big(\|\nabla u\|^2 + \|\nabla\Delta u\|^2\big)\\\nonumber
&&\leq\ C(1 + \mu)(\|u\|_{H^2} + \|V\|_{H^2})\big(\|\nabla u\|^2 +
    \|\nabla\Delta u\|^2 + \|\nabla w\|^2\big).
\end{eqnarray}
Similarly, insert \eqref{i3} and \eqref{i4} into \eqref{g} and
\eqref{h}, and then add the resulting inequalities together, we
obtain that
\begin{eqnarray}\label{fg}
&&\frac{d}{dt}\int_{R^2}\big(|w|^2 + |\Theta|^2\big)dx +
    \mu\big(\|\nabla w\|^2\| + \nabla\Theta\|^2\big)\\\nonumber
&&\leq\ C(\mu^2 + \frac{1}{\mu})(\|u\|_{H^2} +
    \|V\|_{H^2} + \|\Theta\|)\big(\|\nabla u\|^2
    + \|\Delta\nabla u\|^2\\\nonumber
&&\quad +\ \|\nabla w\|^2 + \|\nabla\Theta\|^2) +\
\frac{C}{\mu^3}\|\nabla u\|^2.
\end{eqnarray}
Combining \eqref{ff} with \eqref{fg}, we finally arrive at
\begin{eqnarray}\label{00}
&&\frac{d}{dt}\int_{R^2}\Big[\big(|u|^2
    + |\Delta u|^2 + 2|V|^2 + 2|\Delta V|^2\big) + |w|^2 + |\Theta|^2\Big]dx\\\nonumber
&&\quad +\ \mu\big(\|\nabla u\|^2 + \|\nabla\Delta
    u\|^2 + \|\nabla w\|^2 + \|\nabla\Theta\|^2\big)\\\nonumber
&&\leq\ \frac{C}{\mu^3}\|\nabla u\|^2 + C(\mu^2 +
\frac{1}{\mu})(\|u\|_{H^2} + \|V\|_{H^2} +
    \|\Theta\|)\\\nonumber
&&\quad \times\ \big(\|\nabla u\|^2 + \|\nabla\Delta u\|^2 +
    \|\nabla w\|^2 + \|\nabla\Theta\|^2\big).
\end{eqnarray}

Therefore, if the initial data is sufficiently small, we can find
some positive time $T$, such that
\begin{equation}\label{i5}
\|u\|_{H^2} + \|V\|_{H^2} + \|\Theta\| \leq \frac{\mu^2}{2C(1 +
\mu^3)}
\end{equation}
for all $0 \leq t \leq T$. Moreover, in this case, \eqref{00}
implies that
\begin{eqnarray}\label{i6}
\nonumber&&\frac{d}{dt}\int_{R^2}\Big[\big(|u|^2
    + |\Delta u|^2 + 2|V|^2 + 2|\Delta V|^2\big) + |w|^2 + |\Theta|^2\Big]dx\\\nonumber
&&\quad +\ \mu\big(\|\nabla u\|^2 + \|\nabla\Delta u\|^2 +
\|\nabla w\|^2 + \|\nabla\Theta\|^2\big) \leq
\frac{C}{\mu^3}\|\nabla u\|^2.
\end{eqnarray}
Thus,  noting \eqref{g1}, \eqref{h1}, and by integrating the above
inequality from 0 to $t$, $0 \leq t \leq T$, with respect to $t$,
we obtain
\begin{eqnarray}\label{i7}
&&\big(\|u\|^2_{H^2} + \|V\|^2_{H^2} + \|w\|^2 +
  \|\Theta\|^2\big)\\\nonumber
&&\quad +\ \mu\int_0^\infty\big(\|\nabla u\|_{H^2}^2
    + \|\nabla w\|^2 + \|\nabla\Theta\|^2\big)dt\\\nonumber
&&\leq\ C(\|u_0\|^2_{H^2} + \|V_0\|^2_{H^2}
    + \|w_0\|^2 + \|\Theta_0\|^2)
    + \frac{C}{\mu^3}\int_0^\infty\|\nabla u\|^2\ dt\\\nonumber
&&\leq\ C(1 + \frac{1}{\mu^4})(\|u_0\|^2_{H^2} + \|V_0\|^2_{H^2}
    + \|\nabla\theta_0\|_{H^1}^2),
\end{eqnarray}
where in the last inequality, we used the basic energy law
\eqref{f4}. Then, it follows from \eqref{i7} that if the initial
data satisfies
\begin{equation}\label{i10}
\|u_0\|_{H^2}^2 + \|V_0\|_{H^2}^2 + \|\nabla\theta\|^2_{H^1} \leq
\frac{\mu^8}{M(1 + \mu^{10})}
\end{equation}
for a big enough constant $M$ $(M \geq 2C^3)$ independent of
$\mu$, then \eqref{i5} will still hold for all the latter time $t
\geq T$ with a uniform constant $C$ independent of $t$ and $\mu$.
Furthermore, when \eqref{i10} holds, we can deduce from \eqref{i7}
that
\begin{eqnarray}
&&\big(\|u\|^2_{H^2} + \|V\|^2_{H^2}
  + \|w\|^2 + \|\Theta\|^2\big)\\\nonumber
&&\quad +\ \mu\int_0^\infty\big(\|\nabla u\|_{H^2}^2
    + \|\nabla w\|^2 + \|\nabla\Theta\|^2\big)dt\\\nonumber
&&\leq\ C(\|u_0\|^2_{H^2} + \|V_0\|^2_{H^2}
    + \|w_0\|^2 + \|\Theta_0\|^2)
    + \frac{C}{\mu^3}\int_0^\infty\|\nabla u\|^2\ dt\\\nonumber
&&\leq\ \frac{\mu^4(1 + \mu^4)}{C^2(1 + \mu^{10})}
\end{eqnarray}
provided $C$ is big enough.

At last, we should show \eqref{z1}, \eqref{z2} and \eqref{z3}. By
\eqref{i5} and Sobolev imbedding theorem, we have
$$\|V\|_{L^\infty} \leq c\|V\|_{H^2} \leq \frac{c\mu^2}{2C(\mu^3 + 1)}$$
for appropriate constants $c$ and $C$. Thus, we have
\begin{equation}
\|V\|_{L^\infty} \leq
\begin{cases}
\frac{1}{\mu} \leq 1, \quad {\rm for}\quad \mu \geq 1,\\
1 \leq \frac{1}{\mu}, \quad {\rm for}\quad \mu \leq 1,
\end{cases}
\end{equation}
which means that \eqref{z1} and \eqref{z3} hold. To show
\eqref{z2}, we deduce from \eqref{i5} and \eqref{h1} that
$$\|\nabla\theta\| \leq 2\mu(\|\Theta\| + \|u\|_{H^2})
\leq \frac{\mu^3}{C(1 + \mu^3)} \leq 1.$$

The proof of Theorem 2.2 is completed.

\section*{Acknowledgments} The author would like to thank
professor Chun Liu and professor Yi Zhou for their constructive
suggestions and helpful discussions. The author was partially
supported by the NSFC grant 10801029, the PSFC grants 20070410160
and 200801175.

\bibliographystyle{amsplain}

\end{document}